\documentclass[11pt,letterpaper,reqno]{amsart}
\usepackage[latin1]{inputenc}
\usepackage[T1]{fontenc}
\usepackage{amssymb,amsthm,amsmath,amsfonts,mathrsfs}
\usepackage{hyperref}
\usepackage{graphicx}
\usepackage[english]{babel}
\usepackage{color}

\usepackage{esint}

\usepackage{amsmath}

\usepackage{hyperref}

\usepackage{cite}
\newtheorem{theorem}{Theorem}
\newtheorem{proposition}{Proposition}

\newtheorem{remark}{Remark}[section]
\numberwithin{lemma}{section}
\numberwithin{proposition}{section}
\def\R{\bf \mathbb{R}}
\def\eps{\varepsilon}

\newcommand{\ts}{\textstyle}
\newcommand{\ds}{\displaystyle}

\newcommand{\be}{\begin{equation} \label}
\newcommand{\ee}{\end{equation}}
\newcommand{\N}{\mathbb{N}}
\newcommand{\cal}{\mathcal}

\begin{document}

\tracingpages 1
\title[Optimal Hamilton-type gradient estimates]{Optimal Hamilton-type gradient estimates and large time heat kernel bounds
\\ for noncompact manifolds}
\author[Chabi and Souplet]{Loth Damagui Chabi and Philippe Souplet}
\address{Universit\'e Sorbonne Paris Nord \& CNRS UMR 7539, Laboratoire Analyse G\'eom\'etrie et Applications, 93430 Villetaneuse,  France.}
\email{chabi@math.univ-paris13.fr; souplet@math.univ-paris13.fr }

\begin{abstract}
  We derive localized and global noncompact versions of Ham\-ilton's gradient estimate
 for positive solutions to the heat equation on Riemannian manifolds with Ricci curvature bounded below.
 Our estimates are essentially optimal and significantly improve on all previous estimates of this type.

As a main application, we
obtain a {\it large time} logarithmic gradient estimate for the heat kernel,
which is almost sharp and
considerably improves on previously known results.
Indeed, whereas the precise behavior was known in the small time range,
the large time behavior was rather poorly understood
and remained an essentially open problem, which we here solve to a large extent.
 
 As further applications, we derive a new, space only, local pseudo-Harnack inequality,
as well as estimates of the spatial modulus of continuity of solutions.



\end{abstract}
\maketitle
\section{Introduction}

Logarithmic gradient estimates have played a fundamental role in the study of elliptic and parabolic equations on Riemannian manifolds
since the fundamental works of Cheng-Yau \cite{CY} for 
positive harmonic functions,
 and of Li-Yau \cite{LiYau86} 
and Hamilton~\cite{Ham93} for the homogeneous heat equation in the local and compact cases, respectively.
Some Hamilton-type estimates in the noncompact case were later obtained
in \cite{souplet2006sharp,Kot,DN}.

All three types of estimates for the heat equation (Li-Yau, global and local Hamilton) have been extended in many directions, including various types of linear and nonlinear parabolic equations,
and a number of important applications have been developed, especially heat kernel bounds and Liouville type theorems. 
We refer to, e.g.,~\cite{ATW,BBG,bauer2015li,CCM,colding2021optimal,DN,souplet2006sharp,HZL,HZ,Li91,LX,LZ,Liu,Kot,lin2019ancient,
LNVV,Lu23,MZS,mosconi2021,Wa,WWa,Wu,Ya,ZZ,ZhIMRN,Zh} and to the recent survey paper~\cite{ZhSurvey}
(we stress that this list is necessarily incomplete).
 Let us also mention that the localized Hamilton type estimates have been used in the study of the Ricci flow
(see \cite{Bai, BBG,ZhIMRN}).

In view of the fundamental role of such bounds, it is important to look for the best possible ones. 
However, in spite of this intense activity, it seems that {\it optimal} logarithmic gradient estimates, even for the heat equation, have been missing so far.
 Moreover, regarding logarithmic gradient estimates for the heat kernel,
the precise behavior seems to be known only in the small time range,
while the large time behavior  is rather poorly understood
and remains an essentially open problem.

\smallskip

In this paper, we want to address the following  two questions:
\smallskip
\be{questionKernel}
\begin{aligned}
&\qquad\quad\ \hbox{Can one obtain precise large time}\\
\noalign{\vskip -1mm}
&\hbox{logarithmic gradient estimate for the heat kernel ?}
\end{aligned}
\ee
\be{questionH}
\hbox{What is the best possible localized Hamilton-type gradient estimate ?}
\ee

\smallskip

\noindent Our main results are:

\vskip 1mm
$\bullet${\hskip 1.2mm}an almost complete answer to question~\eqref{questionH},
by establishing a new, nearly optimal Hamilton type gradient estimate 
which covers both the localized and global noncompact cases;

\vskip 1mm

$\bullet${\hskip 1.2mm}as a main application, for manifolds with nonnegative Ricci curvature, a large time logarithmic gradient estimate for the heat kernel,
which is essentially sharp.  
It considerably improves on all previously known results
and solves the open problem \eqref{questionKernel} to a large extent.

\vskip 1mm

Further applications of our Hamilton type gradient estimate, such as spatial pseudo-Harnack inequalities
and estimates of the spatial modulus of continuity of positive solutions will be also derived.

\section{Main results}

\subsection{Logarithmic gradient estimate for the heat kernel}

The manifold $\mathcal{M}$ is said to be complete if every geodesic extends to infinity.
Throughout this paper, 
$G(x, y, t)$ will denote the fundamental solution of the heat equation on $\mathcal{M}$,
the natural logarithm will be denoted by $\log$
and $c(n), C(n)$ will denote positive dimensional constants.

To place our result on the heat kernel in perspective, let us first 
consider the case of compact manifolds, for which the logarithmic gradient estimate
 \be{heatkernelSmallTime00}
\frac{|\nabla_x G(x,y,t)|}{G(x,y,t)}\le C(\mathcal{M})\Big(\frac{d(x,y)}{t}+\frac{1}{\sqrt t}\Big),
 \quad x,y\in \mathcal{M},\ t>0
 \ee
was established in \cite{MaSt96,Hsu,STu} (see also Remark~\ref{JustifD2}).
Comparing with the explicit Gaussian heat kernel, for which
\be{heatkernelRn0}
\frac{|\nabla_x G(x,y,t)|}{G(x,y,t)}= \frac{d(x,y)}{2t},
\ee
we see that \eqref{heatkernelSmallTime00} is far from \eqref{heatkernelRn0} in the limit $d(x,y)/\sqrt{t}\ll 1$. %%%
This difference is quite natural 
since an upper estimate qualitatively similar to \eqref{heatkernelRn0}
would necessarily imply 
\be{nablaG0}
\nabla_xG(y,y,t)=0,
\ee
but \eqref{nablaG0} is not true in general, %%%
and in particular $G(\cdot,y,t)$ need not achieve its maximum at 
 the diagonal point $x=y$ -- see \cite{BM25,BM26} and Proposition~\ref{prop-nonmax}(ii) below
 (although these properties 
 may still hold for some $y$ and special classes of manifolds, see \cite{ACMM,BM25,BM26} and Remark~\ref{JustifTilde}).
However, for given $y,t$, once recentring around a global maximum point $\tilde y=\tilde y(y,t)$ of $G(\cdot,y,t)$, we have the estimate
\be{estimD1}
\frac{|\nabla_x G(x,y,t)|}{G(x,y,t)}\le C(\mathcal{M})\frac{d(x,\tilde y)}{t},\quad x\in \mathcal{M},
\ee
qualitatively similar to \eqref{heatkernelRn0}, as a consequence of the known log-Hessian bound
\be{estimD2}
|D^2_x \log G(x,y,t)| \le C_1(\mathcal{M})\Big(\frac{d^2(x,y)}{t^2}+\frac{1}{t}\Big),\quad x,y\in \mathcal{M},\, t>0
\ee
(cf.~\cite{MaSt96,Hsu,STu} and see also Remark~\ref{JustifD2} for more details).
We note that such log-Hessian bounds play an important role in stochastic analysis on manifolds,
in relation with the study of Brownian bridge and path spaces (see  \cite{CLW,CCT21,NeSa} and the references therein).

For complete noncompact manifolds with nonnegative Ricci curvature,
the best, global logarithmic gradient estimate available is:
 \be{heatkernelSmallTime0}
\frac{|\nabla_x G(x,y,t)|}{G(x,y,t)}\le (1+\delta)\frac{d(x,y)}{2t}+\frac{C(n,\delta)}{\sqrt t},
 \quad x,y\in \mathcal{M},\ t>0
 \ee
 for all $\delta>0$ (see \cite{Kot, LZh} and cf.~also  \cite{STu,Hsu,CLW,NeSa}).
As before, compared with the Euclidean case $\mathcal{M}={\R}^n$,
\eqref{heatkernelSmallTime0} is asymptotic to \eqref{heatkernelRn0},
hence optimal, in the short time limit $d(x,y)\gg \sqrt{t}$,
but far from 
\eqref{heatkernelRn0} in the large time limit $d(x,y)\ll \sqrt{t}$.
 This is again linked to the failure of property \eqref{nablaG0} and to the non maximality of $G$ on the diagonal,
which may still occur in the noncompact case (see Proposition~\ref{prop-nonmax}(i) below), %%%
and we thus need to recenter around a global maximum point. 
 It turns out that, after this recentring, 
 our first main result provides the sought for, almost sharp estimate,
 similar to \eqref{heatkernelRn0} up to a logarithmic factor.

\begin{theorem} \label{propBernAp}
Let $\mathcal{M}$ be a  complete, noncompact manifold of dimension $n\ge 2$ with nonnegative Ricci curvature, 
and let $y\in \mathcal{M},\ t>0$. There exists 
$\tilde y=\tilde y(y,t)\in \mathcal{M}$ such that $\max_{x\in\mathcal{M}}G(x,y,t)=G(\tilde y,y,t)$ and,
for all $\delta>0$, there holds
 \be{heatkernelSmallTimeCor}
 \frac{|\nabla_x G(x,y,t)|}{G(x,y,t)}\le \frac{d(x,\tilde y)}{2t}
 \left\{1+\delta+C(n,\delta)\log\Big(1+\frac{\sqrt{t}}{d(x,\tilde y)}\Big)\right\},\quad x\in \mathcal{M}.
 \ee
\end{theorem}
 
  \begin{remark}
  (i) In view of the inequalities
$$\big|d(x,\tilde y)-d(x,y)\big| \le d(\tilde y, y)\le C(n)\sqrt{t}$$
 (cf.~\eqref{DefTilde}) and $\log(1+X)\le X$, we see that \eqref{heatkernelSmallTimeCor}
 in particular contains~\eqref{heatkernelSmallTime0}.
It improves it considerably in the range $d(x,\tilde y)\ll \sqrt{t}$
(whereas the two inequalities are equivalent and optimal in the range $d(x,\tilde y)\gg \sqrt{t}$).
Theorem~\ref{propBernAp} applies regardless  of whether $\tilde y\ne y$ or $\tilde y=y$, and examples corresponding to either situation
 are given in Proposition~\ref{prop-nonmax} and Remark~\ref{JustifTilde}, respectively.

\smallskip

(ii) Although we do not presently know if \eqref{heatkernelSmallTimeCor} can be improved in general, 
we see 
that the large time estimate, for general noncompact complete manifolds with nonnegative Ricci curvature,
cannot depart from the Euclidean case (where $\tilde y=y$) 
by more than a logarithmic factor in $X=\frac{\sqrt{t}}{d(x,\tilde y)}\gg 1$.

\smallskip

 (iii) Unlike in the compact case, Theorem~\ref{propBernAp} is not a consequence of known log-Hessian estimates.
 Indeed, 
when $\mathcal{M}$ is noncompact, an estimate of the form \eqref{estimD2}
is known to hold only for $x,y,t$ in bounded sets (see \cite{CLW,NeSa,CCT21}),
which cannot yield a uniform estimate such as~\eqref{heatkernelSmallTimeCor}.
A noncompact analogue of \eqref{estimD2} requires very restrictive assumptions on the geometry of $\mathcal{M}$ at infinity 
(see \cite{Ai04} and the introduction of \cite{NeSa}).
  \end{remark}

\subsection{Localized and global noncompact Hamilton-type gradient estimates}
  Theorem~\ref{propBernAp} will be a (nontrivial) consequence of
  our general answer to question \eqref{questionH}.
To give it a precise meaning,  we reformulate \eqref{questionH} as follows. 
Here and throughout this paper, for given $R,T>0$ and $x_0\in \mathcal{M}$,
we denote by $B(x_0,R)$ the geodesic ball of center $x_0$ and radius $R$ on the manifold $\mathcal{M}$,
and by  $Q_{x_0,R,T}$ (or $Q_{R,T}$) the parabolic cylinder $B(x_0,R)\times (0,T]$.
 We will also denote $Q_{\infty,T}=\mathcal{M}\times(0,T]$.
 For given $\mathcal{M}$ and fixed $x_0\in \mathcal{M}$, $R\in(0,\infty]$ and $M,T>0$,
we define $E=E_{M,x_0,R,T}$ by:
$$E=\big\{\hbox{solutions $u$ of the heat equation in $Q_{x_0,R,T}$ with $0<u\le M\big\}$}.$$
Focusing on the ``Hamilton ratio'' $s=\log(M/u)$ as basic variable,
we are looking for the {\it minimal} function $H=H(s,t,R)$ such that, for all $u\in E$,
\be{questionH2}
\frac{|\nabla u|^2}{u^2} \le H\big[\log(M/u),t,R\big] \quad\hbox{in $Q_{R/2,T}$}. 
\ee
In other words, can we determine (up to dimensional multiplicative constants) the quantity
\be{defH}
\begin{aligned}
&H(s,t,R):= \\
&\ \sup \bigg\{\frac{|\nabla u(x,t)|^2}{u^2(x,t)}:\, u\in E, \, x_0\in \mathcal{M},\, x\in B(x_0,R/2),\, \log\big(\ts\frac{M}{u(x,t)}\big)=s\bigg\},
\end{aligned}
\ee
for all $s, t, R>0$~? 
We point out that, even in the global noncompact case $R=\infty$, this question 
does not seem to have been completely answered either
(the previously known results will be recalled in Section~\ref{SecDisc} below).

Our main results  on question \eqref{questionH}, respectively in the localized and global noncompact cases, are the following.

\goodbreak

\begin{theorem} \label{propBern}
Let $\mathcal{M}$ be a complete Riemannian manifold of dimension $n \ge 2$ such that
$Ricci (\mathcal{M}) \ge -k$ for some $k\ge 0$.
For any $R,T>0$, $x_0\in \mathcal{M}$ and any  
positive solution $u$ of the heat equation with
$u \le M$ in $Q_{R,T}$, there holds
 \be{grabound} 
\frac{|\nabla u|^2}{u^2}\le  H_0\Big[\log\Big(\frac{M}{u}\Big),t,R\Big]  
\ \hbox{ in $Q_{R/2,T}$,}
\ee 
with
\be{defgthmBern}
H_0(s,t,R)=
C(n)\begin{cases}
\ds\bigg(\frac{|\log s|^2}{R^2}+\Big(\frac{1}{t}+k\Big)|\log s|\bigg)s^2,&0<s\le 1/2\\
\noalign{\vskip 1mm}
\ds\frac{s^2}{R^2}+\Big(\frac{1}{t}+k\Big)s,&s>1/2.
\end{cases}
\ee
\end{theorem}

\goodbreak
   
 \begin{theorem}\label{Hamilton1}
Let $\mathcal{M}$ be a noncompact complete manifold of dimension $n\ge 2$ 
with $Ricci(\mathcal{M}) \geq -k$ for some $k\geq 0$.
For any  
positive solution $u$ of the heat equation with
$u \le M$ in $\mathcal{M} \times (0, T)$, there holds
\be{Hami}
 \frac{|\nabla u|^2}{u^2} \le H_1\big(\log(M/u)\big),  
\quad x \in \mathcal{M},\ t>0,
\ee
where 
\be{Hami1}
H_1(s)= C(n)\Big(\frac{1}{t} + k \Big) 
\begin{cases}
s^2|\log s|,&0<s\le 1/2\\
\noalign{\vskip 1mm}
s,&s>1/2.
\end{cases}
\ee
 \end{theorem}

\smallskip

The essential optimality of Theorems~\ref{propBern}-\ref{Hamilton1} can be seen from the following lower estimate on $H$,
in the case $\mathcal{M}={\R}^n$,
that we obtain by suitable use of special solutions.

\begin{proposition} \label{prop-optim2}
Let $\mathcal{M}={\R}^n$ and let $H$ be defined by \eqref{defH}.
Then we have 
\be{lowerg}
H(s,t,R) \ge c(n) \bigg(\frac{1}{R^2}+ \frac{1}{t(s+1)}\bigg) s^2,\quad s, t, R>0
\ee
and
\be{lowergB}
H(s,t,\infty) \ge  \frac{c(n)}{t}\frac{s^2}{s+1},\quad s, t>0.
\ee

\end{proposition}

Comparing with the lower bound in \eqref{lowerg}, we see that estimate \eqref{grabound}-\eqref{defgthmBern} in Theorem~\ref{propBern}~is the best possible result of this type up to the $\log\log$ order.
More precisely, it is:

 \vskip 1pt
 
 $\bullet${\hskip 1.2mm}sharp for $s\ge 1/2$
 (in both the spatial and temporal parts)

 \vskip 1pt
 
 $\bullet${\hskip 1.2mm}almost sharp for $s\le 1/2$, up to the order $\log\log(M/u)$.
 \vskip 1pt

 \vskip 1pt

\noindent Detailed comparison with previously known results is given in the next section.

 \smallskip
 
 \subsection{Discussion on Hamilton-type gradient estimates.} 
 \label{SecDisc}
 Theorems~\ref{propBern}-\ref{Hamilton1} significantly improve on all previously known results. 

In the global noncompact 
case ($R=\infty$), under the assumptions of Theorem~\ref{Hamilton1}, 
the bound in \cite{Kot} gave
\be{upperginfty}
H(s,t,\infty)\le\Big(\frac{1}{t} + 2 k \Big) s,
\ee
which was also the form of the celebrated result of Hamilton \cite{Ham93}, given there for compact manifolds without boundary.
While \eqref{upperginfty}  was optimal for $s$ large in view of \eqref{lowergB}, there was a gap between 
\eqref{lowergB}  and \eqref{upperginfty} for $s$ small, which is essentially filled by our Theorem~\ref{Hamilton1}
(up to $\log\log$ order).

Turning to the local case, let $\mathcal{M}$, $u$ be as in Theorem~\ref{propBern}.
We first recall the Li-Yau estimate  \cite{LiYau86}, 
which is of a different, parabolic nature:
  \be{2}
\alpha\frac{|\nabla u|^2}{u^2} -\frac{u_t}{u} \le  C(n,\alpha)
\left(\frac{1}{R^2}+\frac{1}{t}+k\right) 
\ \hbox{ in $Q_{R/2,T}$,}
\ee
with any $\alpha\in(0,1)$.
It has the advantage of being local, but only allows to compare the
solution at different times (it is in particular used in \cite{LiYau86} to derive 
 the parabolic Harnack inequality \eqref{HarnackLY}).
On the contrary, Hamilton's estimate~\eqref{upperginfty} allows to compare the solution at different points at the same time,
but has the drawback of being a global result.
To overcome this, especially for noncompact manifolds, 
a localized estimate of modified Hamilton-type was obtained  in \cite[Theorem~1.1]{souplet2006sharp}
in the form
\be{uppergSZ}
H(s,t,R)\le C(n)\bigg(\frac{1}{R^2}+ \frac{1}{t}+k\bigg)  \big(1+s^2\big).
\ee
Comparing with the lower bound \eqref{lowerg}, the spatial part of \eqref{uppergSZ} (in factor of $R^{-2}$) was sharp for $s$ large,
but there was a gap between \eqref{lowerg} and \eqref{uppergSZ} for $s$ small, as well as for all $s$ in the temporal part (in factor of $t^{-1}$).
Under the extra assumption $m:=\inf_{Q_{R,T}} u>0$, the work \cite[Theorem~1.1]{DN} 
(see Remark~\ref{remDN} below for more details) gave a different estimate of $H$ in terms of the 
additional parameter $\bar s:=\log(M/m)\ge s$:
\be{uppergDD}
H(s,t,R)\le  C(n)\left(\frac{1+\log \bar s}{R^2}+\frac{1}{t}+k\right)s,
\ee
and the same estimate was obtained in \cite[Theorem~1.1]{HuZ} in the case of metric measure spaces
 with Riemannian curvature-dimension condition 
($RCD^*(k,n)$ spaces).
Comparing \eqref{lowerg} and \eqref{uppergDD}, we see that this time the temporal part of \eqref{uppergDD} (in factor of $t^{-1}$)
was sharp for $s$ large, while there was still a gap for $s$ small,
and the spatial parts were not comparable due to the presence of $\bar s$.
Estimate \eqref{grabound}-\eqref{defgthmBern} in our Theorem~\ref{Hamilton1}
fill these gaps (completely for$s$ large and up to $\log\log$ order for $s$ small).
We note in turn that it is stronger than \eqref{uppergSZ} and 
\eqref{uppergDD}, owing to $\log s\le 1+\log^2s$ and $\log s\le\log \bar s$.

To summarize, we see that the optimal behavior  that we obtain for $s$ large answers the natural conjecture that could be made from the partial previous results and that, meanwhile,
the almost optimal $s^2$ behavior for $s$ small, instead of~$s$, represents a strong gain 
which was rather unexpected in view of Hamilton's classical result.

  \begin{remark}
  (i)
The proof of Theorem~\ref{propBern}-\ref{Hamilton1} is based on a refined Bernstein method.
  Different from past studies on gradient bounds via Bernstein techniques, 
we introduce a key new idea to reach our sharp estimates.
Namely, we construct a novel, tricky auxiliary function,
based on nontrivial changes of variables and careful asymptotic analysis,
which lead to an optimal choice in the context of the Bernstein method.
Our proof is organized in such a way as to explain (see Step 3)
the construction of our auxiliary function.

As for the proof of Theorem~\ref{propBernAp}, it combines:

\begin{itemize}
\item[(a)]the first case of \eqref{Hami}-\eqref{Hami1} from Theorem~\ref{Hamilton1}, integrated along a geodesic

\vskip 2pt

\item[(b)]a careful time splitting procedure in which the Harnack inequality
is applied to control the ratio $M/G(\tilde y,y,t)$
where $M$ is the upper bound of $G$ over a suitable past cylinder.
\end{itemize}
We note that, in part (b) of the argument,
using the two-sided Gaussian estimate (as in, e.g.,~\cite{souplet2006sharp})
instead of the Harnack inequality would 
lead to an unaffordable loss of information for large $t$.

\smallskip

(ii) We do not know presently if the logarithmic factors  $|\log\log(M/u)|$ for $u$ close to $M$ in \eqref{grabound}-\eqref{defgthmBern} can be removed
(and similarly for the logarithmic correction in Theorem~\ref{propBernAp}),
or if they are required, at least on some manifolds.
However our estimate is optimal in the scale or powers of $\log(M/u)$, unlike previously known results.

\smallskip

 (iii) Estimate \eqref{grabound}-\eqref{defgthmBern} can be rewritten as a bound for the gradient of a suitable functional. See formulas \eqref{graboundRem} and \eqref{graboundRem2} below.
    \end{remark}

 \subsection{Further applications}

As a consequence of Theorem~\ref{propBern}, we derive the following sharp, space-only, pseudo-Harnack inequality.

\begin{theorem} \label{propBernAppl}
Under the assumptions of Theorem~\ref{propBern}, for each $t\in(0,T]$ and $x,y\in B(x_0,R/2)$, we have
 \be{pseudoHarnack} 
u(y,t)\le L M^{1-\theta}u^\theta(x,t),
\ee
with
 \be{pseudoHarnack2}
 \theta=\big(1+C\ts\frac{d(x,y)}{R}\big)^{-1},
\quad  L=2\exp\big(CR^2(t^{-1}+k)\big),
\ee
where $d$ denotes the geodesic distance and $C=C(n)>0$.
\end{theorem}

The sharpness of \eqref{pseudoHarnack}, at least for $k=0$, can be seen by considering the family of solutions
 $u_\tau(x,t)=G(x_0,x,t+\tau)$ ($\tau>0$), where $G$ is the heat kernel of $\mathcal{M}$
(see Remark~\ref{remOpt1} for details).
Inequality \eqref{pseudoHarnack} improves the analogous (non-sharp) estimate that can be obtained from  \cite[Theorem~1.1]{souplet2006sharp}, 
and was implicitly mentioned in \cite[Remark~2.1]{souplet2006sharp}.
It is to be compared with the classical parabolic Harnack inequality \cite{LiYau86}:
\be{HarnackLY}
u(y,s)\le u(x,t)(t/s)^{n/2} e^{d^2(x, y)/(4(t-s))}, \quad x, y \in \mathcal{M},\  t>s > 0,
\ee
which holds under the same assumptions with $k=0$. 
Here $\theta=1$ but \eqref{HarnackLY} only allows to compare the solution at different times,
contrary to \eqref{pseudoHarnack}.
We stress  that $\theta\to 1$ in~\eqref{pseudoHarnack}-\eqref{pseudoHarnack2} at a ``microscopic'' scale, namely as $d(x,y)/R\to~0$.

\smallskip

As a further consequence of Theorem~\ref{propBern}, we obtain the following uniform estimates 
on the spatial modulus of continuity of positive solutions of the heat equation.
Here we separate the regions far from the supremum $u(x,t)\le M/2$ and close to the supremum
$u(x,t)\ge M/2$.

\begin{theorem} \label{propBernAppl2}
Under the assumptions of Theorem~\ref{propBern}, let $t\in(0,T]$, $x,y\in B(x_0,R/2)$.
Define the modulus of continuity $\psi$ by
$$
\psi(d):=
\left\{\hskip -5mm 
\begin{array}{cccc}
&e^{CdR(\frac{1}{t}+k)}\xi^{\frac{Cd}{R}}-1,\hfill
&\xi:=\ts\frac{M}{u(x,t)},\hfill&\hbox{if $u(x,t)\le \ts\frac{M}{2}$,}\\
\noalign{\vskip 2mm}
&\exp\Big[e^{CdR(\frac{1}{t}+k)}\xi^{\frac{1}{1+Cd/R}} -\xi\Bigr]-1,\hfill
&\xi:=\ts\log\big(\ts\frac{M}{u(x,t)}\big),\hfill&\hbox{if  $u(x,t)\ge \ts\frac{M}{2}$,}\\
\end{array}
\right.$$
which satisfies $\lim_{d\to 0}\psi(d)=0$. We have:
$$
d(x,y)\le \frac{cR}{R^2(t^{-1}+k)+|\log\xi|}
\quad\Longrightarrow\quad
\frac{|u(x,t)-u(y,t)|}{u(x,t)}\le 
\psi\big(d(x,y)\big),$$
where $c, C>0$ are dimensional constants.
\end{theorem}

\begin{remark} \label{remDN}
The estimate in \cite[Theorem~1.1]{DN} was stated under the form
$$
\ts\frac{|\nabla u|^2}{u^2}\le C
\Big\{\frac{\log(M/m)}{R^2}+\frac{1}{t}+k+\frac{k^{1/2}}{R}\Big\} \log(M/u),
$$
and that in \cite[Theorem~1.1]{HuZ} under the same form without the term $\frac{k^{1/2}}{R}$.
However, in the proof of \cite[Theorem~1.1]{DN}, a constant $1$, coming from the term $D^2+1$ in formula (2.12), is missing in formula (2.13) p.5336, 
so that the result is actually 
$$
\ts\frac{|\nabla u|^2}{u^2}\le C
\Big\{\frac{\log(M/m)+1}{R^2}+\frac{1}{t}+k\Big\} \log(M/u).$$
In the proof of \cite[Theorem~1.1]{HuZ}, the additional term $1/R^2$ is missing on p.319
when applying formula (3.10) from Lemma~3.2, since the latter requires $\delta=m/M<1/2$.
 \end{remark}

 \medskip

 We close this section by giving the outline of the rest of the paper. Section~\ref{sect2} is devoted to the proof of Theorem~\ref{propBern}.
Theorem~\ref{Hamilton1} is then derived as a direct consequence.
 Section~\ref{ProofKernel} is devoted to the proof of Theorem~\ref{propBernAp}.
Theorems~\ref{propBernAppl}-\ref{propBernAppl2} are proved in Section~\ref{sect3} and 
Proposition~\ref{prop-optim2} in Section~\ref{sect4}.

\section{Proof of Theorems~\ref{propBern} and \ref{Hamilton1}}\label{sect2}

\subsection{Proof of Theorem~\ref{propBern}}

 The proof is divided in three steps. The first two steps are similar to
previous studies but we give details for convenience and completeness.
Throughout the proof, $C$ will denote a generic positive constant depending only on $n$
and on the parameter $\gamma$ (introduced at the beginning of Step~2).
Let $u$ be as in Theorem~\ref{propBern}
 and assume $x_0=0$  and $M=\sup_{Q_{R,T}} u$ without loss of generality. 

\smallskip

{\bf Step 1.} {\it Gradient equation.}
 Setting $U=\log u$ and $A=\log M=\sup_{Q_{R,T}} U$, 
 we have 
$$U_t-\Delta U=|\nabla U|^2.$$
Let $f$ be a function, to be determined later, such that
\be{fmap}
\begin{aligned}
&\hbox{$f\in C^1([0,\infty))\cap C^3(0,\infty)$,\quad $f'>0$ on $(0,\infty)$,}\\
&\hbox{$f$ maps $[0,\infty)$ onto $[-A,\infty)$.}
\end{aligned}
\ee
Let us put
$$v=f^{-1}(-U),\qquad w=|\nabla v|^2.$$
By direct computation, we see that $v$ satisfies the equation
\begin{equation}  \label{eqLocBernsteinPh}
v_t-\Delta v=\frac{f''}{f'}(v)\,|\nabla v|^2
-{f'}(v)|\nabla v|^2.
\end{equation}
We will derive an equation for $w$. 
First notice that
\begin{align*}
w_t&=2\nabla v\cdot (\nabla v)_t=2\nabla v\cdot \nabla\left( \Delta v+{f''\over f'}(v)\,|\nabla v|^2
-{f'}(v)|\nabla v|^2 \right)\\
&=-2f''\,w^{2}+2\displaystyle\Bigl({f''\over f'}\Bigl)'\,w^2-2\Big(f'-{f''\over f'}\Big)\nabla v\cdot \nabla w+2\nabla v\cdot \nabla(\Delta v),
\end{align*}
 where $\cdot$ denotes the Riemannian inner product.
Therefore, we obtain 
\be{eqcalL}
{\cal L}w=-\Delta |\nabla v|^2+2\nabla v\cdot \nabla (\Delta v)+{\cal N}w \quad\hbox{ in $Q_{R,T}$},
\ee
where 
$$
{\cal L}w:= w_t-\Delta w+b\cdot\nabla w,
\quad b:= 2\Big(f'-{f''\over f'}\Big)\nabla v,
$$
\be{defcalN}
{\cal N}w:=
-2f''w^{2}+2\displaystyle\Bigl({f''\over f'}\Bigl)'\,w^2.
\ee
 Recall the Bochner-Wietzenb\"ock formula: if $q: \mathcal{M}\longrightarrow \R$  is a smooth function, then 
$$
\frac12 \Delta |\nabla q|^2=\nabla q\cdot\nabla (\Delta q)+|D^2q|^2+ Ric_\mathcal{M}( \nabla q,\nabla q),
$$
 where $D^2q$ is the Hessian of $q$. It follows that
$$
 -\Delta |\nabla v|^2 +2\nabla v\cdot \nabla (\Delta v)= -2| D^2 v|^2 - 2  Ric_\mathcal{M}( \nabla v,\nabla v)
 \strut\le - 2| D^2v|^2 +2k|\nabla v|^2,
$$
  and \eqref{eqcalL} yields
 \begin{equation}  \label{eqLocBernsteinPhB}
{\cal L}w\le {\cal N}w +2k w-2|D^2 v|^2
 \quad\hbox{ in $Q_{R,T}$}.
\end{equation}

\smallskip

{\bf Step 2.} {\it Cut-off.}
Let $\gamma\in (0,1)$. Fix any $\tau\in(0,T)$ and let $\eta=\eta(x,t):=\eta_1(t)\tilde{\eta}(x)$, 
where $\tilde \eta \in C^2(\bar{B}_R)$ is a cut-off function with the following properties 
\be{eta0}
0<\tilde \eta\le 1\ \hbox{ in $B_R$},\quad \tilde \eta=1\ \hbox{ in $B_{R/2}$, } 
\quad \tilde \eta=0\ \hbox{ in $B_R\setminus B_{3R/4}$,}
\ee 
\be{eta}
|\nabla \tilde \eta|\le C R^{-1}\tilde \eta^{\gamma},\quad |\Delta \tilde \eta|+\tilde \eta^{-1}|\nabla \tilde \eta|^2\le CR^{-2}\tilde \eta^{\gamma}
\ee
(the existence of $\tilde \eta$ is well known, cf.~\cite{LiYau86}) 
and $\eta_1\in C^1([\tau/4,\tau])$ is a cut-off in time such that:
$$
0\le \eta_1\le 1,\quad \eta_1(\tau/4)=0, \quad \eta_1\equiv 1 \ \hbox{in } [\tau/2,\tau),
\quad |\eta_1'(t)|\le C\tau^{-1}\eta_1^\gamma(t)$$
with $C=C(\gamma)>0$. 

Put 
$$z=\eta w$$
and set $\tilde Q:=B_{3R/4}\times(\tau/4,\tau]$.
In the rest of the proof,
all computations will take place in $\tilde Q\cap\{z>0\}$. 
We have
 \be{eqLocBernsteinPhB2}
{\cal L}z =\eta{\cal L}w+w{\cal L}\eta-2\nabla\eta\cdot\nabla w
\ee
and we also have, in a local orthonormal system:
\be{innereta}
\begin{aligned}
2|\nabla\eta\cdot\nabla w|
&=4\eta_1\big|\sum_i\nabla\tilde\eta\cdot v_i\nabla
v_i\big| \leq 4\eta_1\sum_i \tilde\eta^{-1}|\nabla\tilde\eta|^2 v_i^2+\eta_1\sum_i \tilde\eta|\nabla
v_i|^2 \\
&=4\eta_1\tilde\eta^{-1}|\nabla\tilde\eta|^2 w+\eta|D^2 v|^2.
\end{aligned}
\ee
 Combining  \eqref{eqLocBernsteinPhB}, \eqref{eqLocBernsteinPhB2} and \eqref{innereta}, we obtain
$$ 
{\cal L}z+\eta| D^2 v|^2\leq \eta\, {\cal N}w +2k\eta w+ ({\cal
L}\eta+4\eta_1\tilde\eta^{-1}|\nabla\tilde\eta|^2)w.
$$ 
 Using this and the properties of $\eta$, we deduce that
$$
|(b\cdot\nabla\eta)w|
\leq CR^{-1}\eta^\gamma\Bigl|f'-{f''\over f'}\Bigr|w^{3/2}
$$
and then
\begin{equation}  \label{eqLocBernsteinPh2}
{\cal L}z+\eta|D^2v|^2\leq \eta\, {\cal N}w +
(2k\eta+C\tau^{-1}\eta^\gamma +CR^{-2}\eta^\gamma) w+CR^{-1}\eta^\gamma\Bigl|f'-{f''\over
f'}\Bigr|w^{3/2}.
\end{equation}

{\bf Step 3.} {\it Choice of auxiliary function and completion of proof.}
 Without loss of generality, we may look for $f$ in the following, convenient way.
Let $g\in C([0,\infty))\cap  C^2(0,\infty)$, to be determined below, satisfy
\be{hypg1}
 g>0\quad \hbox{ in $(0,\infty)$}
\ee
and
\be{hypg2}
\int_0^1\frac{d\tau}{g(\tau)}<\infty,\quad \int_1^\infty\frac{d\tau}{g(\tau)}=\infty.
\ee
We set
\be{choiceh}
h(s)=\int_0^s\frac{d\tau}{g(\tau)}
\ee
and then
\be{choicef}
f(s)=h^{-1}(s)-A,\quad \hbox{ hence}\ v=h(A-U).
\ee
Note that $v\geq 0$ in $Q_{R,T}$ due to $u\leq M$ and $g>0$ in $[0,\infty)$. 
With this
choice, we see that \eqref{fmap} is satisfied, and we have
\begin{align*}
f'&=g\circ h^{-1},\qquad \quad
f''=\big(g'g\big)\circ h^{-1},\\
\frac{f''}{f'}&=g'\circ h^{-1},\qquad \quad
\left(\frac{f''}{f'}\right)'= \big(gg''\big)\circ h^{-1}.
\end{align*}
Inequality \eqref{eqLocBernsteinPh2} then implies
\be{ineqcalLz1}
{\cal L}z\leq \eta\, {\cal N}w +
\underbrace{\big(C(\tau^{-1}+R^{-2})\eta^\gamma +2k\eta\big)w}_{T_1}
+\underbrace{CR^{-1}\eta^\gamma w^{3/2}|g-g'|}_{T_2}, 
\ee
whereas
\be{ineqcalLz2}
{\cal N}w=2 g(g''-g')w^2. 
\ee
 Here we have omitted without risk of confusion the variables $h^{-1}(v)=A-U$ in the terms 
 $|g-g'|$ and $g(g''-g')$.
We note that the change of variables \eqref{choiceh}-\eqref{choicef}
amounts to taking as new unknown $g=\frac{1}{((A+f)^{-1})'}$ instead of $f$,
with the advantage of replacing the expressions involving $f,f',f''$ in \eqref{defcalN}, \eqref{eqLocBernsteinPh2}
by more transparent and tractable quantities.

 Our aim is now to control each of the positive terms $T_1, T_2$ by $\ts\frac12\eta{\cal N}w$,
which will be designed as to be negative.
Since we restrict to the set $\{z>0\}$ (where $0<\eta\le 1$),  \eqref{ineqcalLz1}-\eqref{ineqcalLz2} may be rewritten as
\begin{align*}
{\cal L}z
\le  2 \eta^{-1}g(g''-g') z^2
+\big(C(\tau^{-1}+R^{-2})\eta^{\gamma-1} +2k\big) z
+CR^{-1}\eta^{\gamma-\frac{3}{2}}  |g-g'| z^{{\frac{3}{2}}},
\end{align*}
hence
\be{eqLzN}
{\cal L}z\le \eta^{-1}z^2\big(\tilde T_1+ \tilde T_2\big),
\ee
where
$$ 
\begin{aligned} 
\tilde T_1&=g(g''-g')+\big(C(\tau^{-1}+R^{-2})\eta^{\gamma} +2k\eta \big)z^{-1},\\
\tilde T_2&=g(g''-g')+CR^{-1}\eta^{\gamma-\frac{1}{2}}  |g-g'| z^{-{\frac{1}{2}}}.
\end{aligned}
$$
Upon taking $\gamma=\frac12$, we have
\be{eqLzN2}
\left\{\
\begin{aligned}
\tilde T_1
&\le g(g''-g')+C\big(\tau^{-1}+R^{-2}+k \big)z^{-1}\\
\tilde T_2
&\le g(g''-g')+CR^{-1} |g-g'| z^{-{\frac{1}{2}}}.
\end{aligned}
\right.
\ee
In view of ensuring the negativity of $\tilde T_1, \tilde T_2$ at points where $z$ is large, 
we thus look for a function $g$ satisfying 
\be{choiceg0}
 F_1\ge CF_2,
\ee
 where
$$F_1=g(g'- g''),\quad F_2=R^{-1}|g-g'|+B,\quad B:=R^{-2}+\tau^{-1}+k.$$
We write $f_1\approx f_2$ if there exist constants $c_1, c_2>0$ (depending only on $n$) such that $c_1 f_1\le f_2\le c_2 f_1$ on $(0,\infty)$
(or on a specified range).

We claim that inequality \eqref{choiceg0} can be {\it optimally} solved, i.e.~$F_1\approx F_2$, provided we make the precise choices
$$g(s)=R^{-1}g_1+Kg_2,$$
with
$$g_1= s\log(e+\ts\frac{1}{s}),\quad g_2=\frac{s}{(s+1)^{1/2}}\log^{1/2}(e+\frac{1}{s})$$
and
\be{choiceKnew}
K=\tau^{-1/2}+k^{1/2}.
\ee
 By direct computations, on $(0,\infty)$, we obtain
$$g_1'\approx\ts \log(e+\ts\frac{1}{s}),\quad
g_1''\approx 
-\frac{1}{s(s+1)^2},$$
$$g_2'\approx\ts\frac{1}{(s+1)^{1/2}}\log^{1/2}(e+\frac{1}{s}),
\quad g_2''\approx\ts\frac{-1}{s(s+1)^{1/2}\log^{1/2}(e+\frac{1}{s})}.$$
Considering the ranges $s\ge 1$ and $0<s\le 1$ separately, we see that:

\vskip 2pt

\noindent $\bullet$ For $s\ge 1$,
$$\begin{aligned}
F_1&\approx (R^{-1}s+K s^{1/2})(R^{-1}+K s^{-1/2})\\
&\approx R^{-2}s+KR^{-1}s^{1/2}+K^2,\\
F_2&\approx R^{-1}\big(R^{-1}s+Ks^{1/2}\big)+R^{-2}+\tau^{-1}+k\\
&\approx R^{-2}s+KR^{-1}s^{1/2}+\tau^{-1}+k.
\end{aligned}$$
\noindent $\bullet$ For $s\le 1$,
$$
\begin{aligned}
F_1
&\approx \Big(R^{-1}s\log(e+\ts\frac{1}{s})+Ks\log^{1/2}(e+\ts\frac{1}{s})\Big)\Big(R^{-1}s^{-1}+K\frac{s^{-1}}{\log^{1/2}(e+\frac{1}{s})}\Big) \\
&=\Big(R^{-1}\log(e+\ts\frac{1}{s})+K\log^{1/2}(e+\ts\frac{1}{s})\Big)\Big(R^{-1}+\frac{K}{\log^{1/2}(e+\frac{1}{s})}\Big) \\
&\approx R^{-2}\log(e+\ts\frac{1}{s})+KR^{-1}\log^{1/2}(e+\ts\frac{1}{s})+K^2,\\
F_2
&\approx R^{-1}\Big(R^{-1}\log(e+\ts\frac{1}{s})+K\log^{1/2}(e+\ts\frac{1}{s})\Big)+R^{-2}+\tau^{-1}+k\\
&\approx R^{-2}\log(e+\ts\frac{1}{s})+KR^{-1}\log^{1/2}(e+\ts\frac{1}{s})+\tau^{-1}+k.
\end{aligned}
$$
In view of \eqref{choiceKnew}, we thus have $F_1\approx F_2$ on $(0,\infty)$, which proves our claim.

It follows  from \eqref{eqLzN}-\eqref{choiceg0} that there exists $L=L(n)>0$ 
such that $\tilde{T}_i\strut< 0$ 
whenever $z\ge L$, hence
\be{PM}
{\cal L} z \strut < 0
\quad\hbox{ in $Q'_L:=\{(x,t)\in \tilde Q;\ z(x,t)\geq L\}$}.
\ee
 Since $z=0$ on the parabolic boundary of $\tilde Q$, we may
 suppose the  maximum of $z$ is reached at $(x_1, t_1)\in Q'_L$. 
 By \cite{LiYau86}, we can assume, without loss of generality,  that $x_1$ is not in the cut-locus of $\mathcal{M}$.  
Then at this point, one has
 $\Delta z \le 0$, 
 $z_t \ge 0$ and $\nabla z=0$, and using \eqref{PM}, we obtain
 $$
  0\le z_t-\Delta z+b\cdot \nabla z< 0:
 $$
a contradiction.
Consequently, $z\le L$ in $\tilde Q$. 
Finally, using $z=\eta|\nabla v|^2$, \eqref{choiceh}, \eqref{choicef} and the definition of $\eta$, this implies
$$|\nabla U|=|\nabla v|\, g(A-U)\le L^{1/2} g(A-U)$$
in $B_{R/2}\times[\tau/2,\tau]$.
 Applying this for each $t\in (0,T)$ with $\tau=t$ 
  and separating the ranges $s\le 1/2$ and $s>1/2$ yields \eqref{grabound}-\eqref{defgthmBern}.
\qed

\subsection{Proof of Theorem~\ref{Hamilton1}}

Let $x_0\in \mathcal{M}$ and $0<t_0\le T$. 
Fixing $R>0$ and applying Theorem~\ref{propBern} on $Q:=B(x_0,R)\times(0,T]$, we obtain
$$\frac{|\nabla u(x_0,t)|}{u(x_0,t)}\le H_0\big(\ts\log\big(\frac{M}{u(x_0,t)}\big),t,R\big)$$
where $H_0$ is defined in \eqref{defgthmBern}.
Since $\mathcal{M}$ is non compact, we may let $R\to \infty$ 
and we get 
$$\frac{|\nabla u(x_0,t)|}{u(x_0,t)}\le \lim_{R\to\infty} H_0\big(\ts\log\big(\frac{M}{u(x_0,t)}\big),t,R\big)=
C\big(\frac{1}{t} + k \big)H_1\big(\ts\log\big(\frac{M}{u(x_0,t)}\big)\big),$$
which concludes the proof. \qed

\section{Proof of Theorem~\ref{propBernAp}}

\label{ProofKernel}

Before giving the proof of Theorem~\ref{propBernAp}, we begin with the following proposition
which shows that property \eqref{nablaG0} is not true in general.

\begin{proposition} 
\label{prop-nonmax}
Let $\mathcal{M}$ be the complete noncompact manifold given by the Riemannian product $\mathcal{M}= \mathcal{M}_0\times {\R}^d$,
where $d\in \N^*$ and $(\mathcal{M}_0,g)$ is a compact connected Riemannian manifold.

\begin{itemize}
\item[(i)] If the first nonzero eigenvalue of $-\Delta_g$ is simple, then
there exist $y\in \mathcal{M}$ and $t_0>0$ such that
$\nabla_x G(y,y,t)\ne 0$ for all $t\ge t_0$.

\vskip 2pt 
\item[(ii)] Assertion (i) remains valid for $\mathcal{M}= \mathcal{M}_0$.
\end{itemize}
\end{proposition}

Recall that the spectral simplicity assumption in Proposition~\ref{prop-nonmax} is generically true by \cite{Uhl},
and it is for instance true for suitable ellipsoids $\mathcal{M}_0$ \cite{EK22}.
The proof of Proposition~\ref{prop-nonmax} is based on modifications of the ideas from \cite[Theorem~2.2]{BM25},
which showed in the case $\mathcal{M}$ compact that the maximum of $G$ is not taken on the diagonal in general.

\begin{proof} 
(i) Write $x=(\bar x,x')\in \mathcal{M}_0\times {\R}^d$.
The heat kernel of $\mathcal{M}$ can be written as
$$G(x,y,t)= K(\bar x,\bar y,t)G_0(x',y',t)$$
where $G_0(x',y',t)=(4\pi t)^{-d/2}e^{-|x'-y'|^2/4t}$
and $K$ is the heat kernel of $\mathcal{M}_0$.
The latter can be written as the spectral expansion 
$$
K(\bar x,\bar y, t)=
\sum_{j=0}^\infty e^{-\lambda_j t} \phi_j(\bar x)\phi_j(\bar y),
$$
where $(\phi_j)_{j\ge 0}$ is an orthonormal sequence of eigenfunctions in $L^2(\mathcal{M}_0,d{\rm Vol}_g)$,
with corresponding eigenvalues $\lambda_0=0<\lambda_1<\lambda_2\le\lambda_3\le\dots$.
Here the eigenvalues are repeated according to multiplicity, and $\lambda_1<\lambda_2$ by assumption.
As a consequence of Sobolev inequalities and Weyl's asymptotic law for the eigenvalues,
the series as well as its gradient converges absolutely pointwise
(see \cite[Section~5]{BM25} for more details and references).
Note that the eigenfunction $\phi_0$ associated with $\lambda_0=0$ is constant,
and $\phi_1$  is non-constant. 
Therefore, there exists $\bar y\in  \mathcal{M}_0$ such that $\nabla(\phi^2_1(\bar y))=2\phi_1(\bar y)\nabla\phi_1(\bar y)\ne 0$.
Noticing that, for all $t>1$,
$$\begin{aligned}
\Big|e^{\lambda_1 t}\nabla_{\bar x} K(\bar y,\bar y,t)&-\phi_1(\bar y)\nabla\phi_1(\bar y)\Big|
=\Big|\sum_{j=2}^\infty e^{-(\lambda_j-\lambda_1) t} \phi_j(\bar y)\nabla\phi_j(\bar y)\Big|\\ 
&\le  e^{-(\lambda_2-\lambda_1)(t-1)} \sum_{j=2}^\infty e^{-(\lambda_j-\lambda_1)} |\phi_j(\bar y)\nabla\phi_j(\bar y)|, 
\end{aligned}$$
we deduce that $e^{\lambda_1 t}\nabla_{\bar x} K(\bar y,\bar y,t)\to\phi_1(\bar y)\nabla\phi_1(\bar y)\ne 0$ as $t\to\infty$.
Consequently, setting $y=(\bar y,0)$ and using $\nabla_{x'}G_0(0,0,t)=0$, there exists $t_0>1$ such that 
 $\nabla_x G(y,y,t)=G_0(0,0,t)\nabla_{\bar x} K(\bar y,\bar y,t)\ne 0$ for all $t\ge t_0$,
and the conclusion follows.

\smallskip

(ii) This follows immediately from the above proof since $G=K$.
\end{proof}

  \begin{remark} \label{JustifTilde}
  Beside the Euclidean case, the equality $\tilde y=y$ (hence in particular  \eqref{nablaG0})
  nevertheless %%%%
  holds under suitable symmetry assumptions on $\mathcal{M}$ and for certain points $y$.
    This is the case, for instance, if $\mathcal{M}$ is a smooth connected hypersurface of revolution around the $x_n$ axis without boundary,
  and $y$ is a point of intersection of $\mathcal{M}$ and the $x_n$ axis
  (cf.~\cite[Theorem~2 and Remark~1]{ACMM} and the references therein, and see also \cite{BM25} for additional examples).
  \end{remark}

\begin{proof}[Proof of Theorem~\ref{propBernAp}]
First note that the claimed existence of the point
$\tilde y=\tilde y(y,t)\in \mathcal{M}$ follows from the two-sided Gaussian bound on $G$ \cite{LiYau86}:
\be{group2}
\frac{c_1}{|B(y, \sqrt{t})|} e^{- d(x, y)^2/c_1t} \le G(x,y, t) \le \frac{c_2}{|B(y, \sqrt{t})|} e^{- d(x, y)^2/c_2t},
\ee
 for all $x, y \in \mathcal{M}$ and $t>0$, with $c_i=c_i(n)>0$.
Indeed, if $d(x, y)>c_3\sqrt{t}$ with $c_3:=(c_2\log(c_2/c_1))^{1/2}$, then
$$G(y,y,t)\ge \frac{c_1}{|B(y, \sqrt{t})|}  \ge \frac{c_1}{c_2}e^{d(x, y)^2/c_2 t}G(x,y, t) >G(x,y, t).$$
Consequently there exists a point $\tilde y(t)=\tilde y(t,y)$ where 
the continuous function $G(\cdot,y,t)$ attains a global maximum, and moreover there holds
\be{DefTilde}
d(y,\tilde y)\le c_3(n)\sqrt{t}.
\ee 

\smallskip

Next, in view of the proof of \eqref{heatkernelSmallTimeCor}, we fix $x, y\in \mathcal{M}$ and $t>0$. 
Set $D=d(x,\tilde y(t))$. It suffices to prove that
\be{hyptD2A2}
\frac{|\nabla_x G(x,y,t)|}{G(x,y,t)}\le C(n)Dt^{-1}\log(\sqrt{t}D^{-1})
\quad\hbox{for all $t\ge AD^2$,}
\ee
where $A=A(n)\ge e$ will be chosen below.
Indeed, if $0<t\le AD^2$ then we have
$d(x,y)\le D+d(\tilde y(t),y)\le (1+c_3\frac{\sqrt{t}}{D})D$ by \eqref{DefTilde}. 
Also, there exists $C_1(A)>0$ such that $C_1(A)\log(1+X)\ge X$ for all $X\in (0,A]$.
Therefore,  for all $\delta\in (0,1]$, \eqref{heatkernelSmallTime0} implies
$$\begin{aligned}
\frac{|\nabla_x G(x,y,t)|}{G(x,y,t)}
&\le (1+\delta)\frac{d(x,y)}{2t}+\frac{C(n,\delta)}{\sqrt t} 
\le \frac{D}{2t}\Big(1+\delta+2(c_3+C(n,\delta))\frac{\sqrt{t}}{D}\Big)\\
&\le \frac{D}{2t}\Big\{1+\delta+2C_1(A)(c_3+C(n,\delta))\log\Big(1+\frac{\sqrt{t}}{D}\Big)\Big\}.
\end{aligned}$$
Let $\gamma~:~[0,1]\to \mathcal{M}$ be a geodesic such that $\gamma(0)=\tilde y(t)$ and $\gamma(1)=x$.
Set $\theta=1/n$.
We introduce a parameter $\eps\in(0,1/4]$ (depending on $x,y,t$) which will be selected below,
and set $t_0=(1-\theta\eps)t$ and $T=\theta\eps t$.

Set $Q=\mathcal{M}\times[t_0,\infty)$.
Since $m(s)=G(\tilde y(s),y,s)$ is nonincreasing in $s$ by the maximum principle, 
we deduce that
$$M:=\sup_{(z,s)\in Q}G(z,y,s)=G(\tilde y(t_0),y,t_0).$$
We claim that
\be{MGxxt}
\log\Big(\frac{M}{G(\tilde y(t),y,t)}\Big)\le\eps.
\ee
By the Harnack inequality \eqref{HarnackLY} from \cite{LiYau86}, we have 
$$
\frac{M}{G(\tilde y(t),y,t)}\le \frac{M}{G(\tilde y(t_0),y,t)}
=\frac{G(\tilde y(t_0),y,t_0)}{G(\tilde y(t_0),y,t)}\le (t/t_0)^{n/2}=(1-\theta\eps)^{-n/2}. 
$$
Consequently, using $-\log(1-X)\le 2X$ 
in $[0,1/2]$, we obtain
$$\log\Big(\frac{M}{G(\tilde y(t),y,t)}\Big)\le-(n/2)\log(1-\theta\eps)\le n\theta\eps\le\eps,$$
which proves \eqref{MGxxt}.

We shall next combine \eqref{MGxxt} with Theorem~\ref{Hamilton1} to control 
$\log\big(\frac{M}{G(x,y,t)}\big)$ by $2\eps$ for suitable choice of $\eps$.
Let 
\be{defsigma0}
\sigma_0=\max\Big\{\sigma\in(0,1],\ \sup_{\tau\in[0,\sigma]}\log\Big(\frac{M}{G(\gamma(\tau),y,t)}\Big)\le 2\eps\Big\}.
\ee
We have $\sigma_0\in(0,1]$ by continuity. Set $\Sigma=\gamma([0,\sigma_0])$ and
$$\bar x=\gamma(\sigma_0),\quad
u(z,s):=G(z,y,t_0+s),\quad v(z)=\log(M/u(z,T)).$$
Then $u(z,T)\ge M/2$ and $v(z)\le 1/2$ for all $z\in\Sigma$.
Applying Theorem~\ref{Hamilton1} 
to the function $u$, which is a positive solution 
of the heat equation in $\mathcal{M}\times[0,T]$, we obtain
\be{nablavz}
|\nabla v(z)|\le C(n)T^{-1/2}v(z)|\log v(z)|^{1/2},\quad z\in\Sigma
\ee
hence
$$\big|\nabla\big(|\log v|^{1/2}\big)\big|\le C(n)(\eps t)^{-1/2},\quad z\in\Sigma.$$
Integrating this along the part $\Sigma$ of the geodesic, 
and since $v(\tilde y(t))\le \eps$ by \eqref{MGxxt}, we deduce that
$$\begin{aligned}
|\log v(\bar x)|^{1/2}
&\ge |\log v(\tilde y(t))|^{1/2}-C(n)(\eps t)^{-1/2}d(\tilde y(t),\bar x)\\
&\ge |\log\eps|^{1/2}-C(n)(\eps t)^{-1/2}D,
\end{aligned}$$
hence
$$|\log v(\bar x)|
\ge |\log\eps|-2C(n)\big(\eps^{-1}|\log\eps|D^2t^{-1}\big)^{1/2}.$$
Choose
$$\eps=C_1D^2t^{-1}\log(tD^{-2})$$
with $C_1=C_1(n)=1+16\,C^2(n)$.
There exists $A=A(n)\ge e$ such that 
$t\ge AD^2$ implies $\eps\le 1/4$.
Since $D^2t^{-1}\le\eps<1$, we then have
\be{logeps}
|\log\eps|\le \log(tD^{-2}),
\ee
hence
$\eps^{-1}|\log\eps|D^2t^{-1}\le 1/C_1$, so that
$$\begin{aligned}
-\log v(\bar x)
&\ge -\log\eps-2C(n)\big(\eps^{-1}|\log\eps|D^2t^{-1}\big)^{1/2}
\ge -\log\eps-\frac12.
\end{aligned}$$
Consequently, $v(\bar x)\le \eps\sqrt{e}<2\eps$, 
and we deduce from \eqref{defsigma0} that $\sigma_0=1$, i.e.~$\bar x=x$, hence $v(x)\le 2\eps$.

To conclude, we go back to \eqref{nablavz} with $z=x$.
Using also \eqref{logeps} and noting that the function $h(X)=X|\log X|^{1/2}$ is increasing on $(0,1/2]$, 
we deduce that
$$\begin{aligned}
\frac{|\nabla_x G(x,y,t)|}{G(x,y,t)}
&=|\nabla v(x)|\le 2C(n)\big(\eps t^{-1}|\log (2\eps)|\big)^{1/2}\\
&\le 2C(n)\big(C_1D^2t^{-2}\log(tD^{-2})\big)^{1/2} \big(\log(tD^{-2})\big)^{1/2}\\
&\le C_2(n) Dt^{-1}\log(tD^{-2})=2C_2(n)Dt^{-1}\log(\sqrt{t}D^{-1}),
\end{aligned}$$
hence \eqref{hyptD2A2}. 
This completes the proof of Theorem~\ref{propBernAp}.
  \end{proof}

    \begin{remark} \label{JustifD2}
    Let $\mathcal{M}$ be a compact manifold.
  First note that \eqref{heatkernelSmallTime00} and \eqref{estimD2}
  are stated in \cite{MaSt96,Hsu,STu} for $t\le 1$ but are easily seen to remain true for all $t>1$ as a consequence of the spectral expansion of the heat kernel
(and the estimates actually have an exponential decay in time for $t\to\infty$).

 Let us next give the simple 
 derivation of \eqref{estimD1} from \eqref{estimD2}.
Here we set $R={\rm diam}(\mathcal{M})$ and  denote by $C$ a generic positive constant depending only on~$\mathcal{M}$.
 By \cite[Theorem 5.4.12 and Section 5.6.3]{LSC02}, the two-sided Gaussian bound \eqref{group2}, hence estimate \eqref{DefTilde},
 are still valid for all $t\in(0,1]$
 (with $c_1,c_2>0$ depending on $\mathcal{M}$ and $c_3(n)$ replaced by $C$).
Since $d(\tilde y,y)\le R$, we also have~\eqref{DefTilde} for $t>1$, hence for all $t>0$.
Fix $t>0$ and $y\in \mathcal{M}$. 
We deduce from \eqref{estimD2} and \eqref{DefTilde} that, for all $z\in B_{R\sqrt{t}}(\tilde y)$,
$$ 
|D^2_x \log G(z,y,t)| 
\le\frac{C}{t}\Big(1+\frac{d^2(z,y)}{t}\Big)
\le\frac{C}{t}\Big(1+\frac{d^2(z,\tilde y)+d^2(\tilde y,y)}{t}\Big)
\le\frac{C}{t}.
$$
 If $d(x,\tilde y)<R\sqrt{t}$, then \eqref{estimD1} follows by using $\nabla_x G(\tilde y,y,t)=0$ and
 integrating this inequality along a geodesic connecting $x$ and $\tilde y$.
 If $d(x,\tilde y)\ge R\sqrt{t}$,  then  \eqref{estimD1}  directly follows from \eqref{heatkernelSmallTime00} and \eqref{DefTilde}:
$$\frac{|\nabla_x G(x,y,t)|}{G(x,y,t)}\le C\Big(\frac{d(x,y)}{t}+\frac{1}{\sqrt t}\Big)
\le C\Big(\frac{d(x,\tilde y)+d(\tilde y,y)}{t}+\frac{1}{\sqrt t}\Big)
\le C\frac{d(x,\tilde y)}{t}.$$
  \end{remark}

\section{Proof of Theorems~\ref{propBernAppl}-\ref{propBernAppl2}}\label{sect3}

We start with the following:

 \begin{remark} 
By direct calculation we see that estimates \eqref{grabound}-\eqref{defgthmBern} can 
be rewritten as:\footnote{ The threshold between the two regions of \eqref{grabound} is $u=M/\sqrt{e}$ but it is easy to see that it can be replaced by $M/2$ for instance.}
\be{graboundRem}
CR^{-1}\ge
\begin{cases}
\Big|\nabla\Big[\log\big(\sqrt{\log(M/u)}+K\big)\Big]\Big|,&\hbox{if $u(x,t)\le M/2$},\\ 
\noalign{\vskip 1mm}
\Big|\nabla\Big[\log\big(\sqrt{|\log\log(M/u)|}+K\big)\Big]\Big|,&\hbox{if $u(x,t)\ge M/2$}, 
\end{cases}
\ee
where $K=R\sqrt{t^{-1}+k}$,
or, alternatively, under the global form
\be{graboundRem2}
\Big|\nabla\Big[\log\Big(\sqrt{\log(M/u)+\log\big(1+\ts\frac{1}{\log(M/u)}\big)}+K\Big)\Big]\Big|\le CR^{-1}.
\ee
 \end{remark} 
 
 \subsection{Proof of Theorem~\ref{propBernAppl2}}

In this proof $C$ denotes a generic positive constant depending only on $n$.
We first note that, for all $a,b,\mu\ge 0$ and $\lambda\ge 1$, 
\be{implelem}
a+\mu\le \lambda(b+\mu) \Longrightarrow a^2\le \lambda^4 b^2+(\lambda^2-1)\mu^2.
\ee
Indeed, we may assume $a\ge\lambda^2b$ (since otherwise it is immediate), 
and the assumption in~\eqref{implelem} then implies
$$a^2+\mu^2=(a+\mu)^2-\lambda^2(b+\mu)^2+\lambda^2(b^2+\mu^2)+2\mu(\lambda^2b-a)\le \lambda^2(b^2+\mu^2),$$
hence $a^2\le \lambda^2 b^2+(\lambda^2-1)\mu^2\le \lambda^4 b^2+(\lambda^2-1)\mu^2$.

Fix $x\in B(x_0,R/2)$. Assume $u(x,t)\le M/2$ (resp., $u(x,t)\ge M/2$) and let
\be{defrhoH}
\rho=\sup\Big\{r\in(0,R/2];\,u(\cdot,t)\le \ts\frac{3M}{4} \ \hbox{(resp., $u(\cdot,t)\ge \frac{M}{e}$)\ in $\overline B(x,r)$}\Big\}.
\ee
By continuity we have $\rho>0$. Pick $y\in B(x_0,R/2)$ with $y\ne x$ and $d(x,y)\le \rho$.
Let $\gamma$ be a geodesic connecting $x$ and $y$ and note that $d(x,z)\le\rho$ for all $z\in\gamma$.

First consider the case $u(x,t)\le M/2$ and
set either $(x_1,x_2):=(x,y)$ or $(x_1,x_2):=(y,x)$ in such a way that $u(x_1,t)\ge u(x_2,t)$.
Set
\be{notationKlambda}
K=R\sqrt{t^{-1}+k},\quad d=d(x,y),\quad \lambda=e^{Cd/R}>1.
\ee
Observing that $M/2$ can be replaced by $3M/4$ in the first part of \eqref{graboundRem}
and integrating along $\gamma$, we obtain
$$\log\Big(\sqrt{\log\big(\ts\frac{M}{u(x_2,t)}\big)}+K\Big)\le 
\log\Big(\sqrt{\log\big(\ts\frac{M}{u(x_1,t)}\big)}+K\Big)+C\ts\frac{d}{R},$$
hence
$$\sqrt{\log\big(\ts\frac{M}{u(x_2,t)}\big)}+K\le \lambda\Big(\sqrt{\log\big(\ts\frac{M}{u(x_1,t)}\big)}+K\Big).$$
Using \eqref{implelem}, we deduce
\be{modulus1a}
\log\big(\ts\frac{M}{u(x_2,t)}\big)\le \lambda^4 \log\big(\ts\frac{M}{u(x_1,t)}\big)+(\lambda^2-1)K^2.
\ee
We note for further reference that \eqref{modulus1a} is true whenever $u(\cdot,t)\le \ts\frac{3M}{4}$
along a geodesic connecting $x_1$ and $x_2$.
Next using $d\le R$, we get
$$\begin{aligned}
\log\big(\ts\frac{u(x_1,t)}{u(x_2,t)}\big)
&=\log\big(\ts\frac{M}{u(x_2,t)}\big)-\log\big(\ts\frac{M}{u(x_1,t)}\big)\\
&\le \lambda^4 \log\big(\ts\frac{M}{u(x_1,t)}\big)+(\lambda^2-1)K^2-\log\big(\ts\frac{M}{u(x_1,t)}\big)\\
&\le C\ts\frac{d}{R} \Big(\log\big(\ts\frac{M}{u(x_1,t)}\big)+R^2\big(t^{-1}+k\big)\Big).
\end{aligned}$$
It follows that
\be{modulus1}
\ts\frac{|u(y,t)-u(x,t)|}{u(x,t)}\le \ts\frac{u(x_1,t)}{u(x_2,t)}-1\le e^{CdR(t^{-1}+k)}\big(\ts\frac{M}{u(x,t)}\big)^{\frac{Cd}{R}}-1.
\ee
Now, if we assume
$$d(x,y)\le \frac{cR}{R^2(t^{-1}+k)+\log(\frac{M}{u(x,t)})}$$
with $c=c(n)>0$ small, we in particular deduce from \eqref{modulus1} that $\frac{u(y,t)}{u(x, t)}\le 4/3$.
Consequently, $\rho=R/2$ in \eqref{defrhoH}
and this concludes the proof in the first case.

Next consider the case $u(x,t)\ge M/2$.
Observe that $M/2$ can be replaced by $M/e$ in the second part of \eqref{graboundRem}.
Also, we may replace $M$ by $N=M+\eps$, so that $u<N$ in $Q_{R,T}$.
Set either 
$(x_1,x_2):=(x,y)$ or  $(x_1,x_2):=(y,x)$ in such a way that 
 $N/e\le u(x_1,t)\le u(x_2,t)$, hence 
\be{monotxi}
0<\log\big(\ts\frac{N}{u(x_2,t)}\big)\le \log\big(\ts\frac{N}{u(x_1,t)}\big)\le 1
\ee
 and
$-\log\log\big(\ts\frac{N}{u(x_2,t)}\big)\ge -\log\log\big(\ts\frac{N}{u(x_1,t)}\big)\ge 0$.
Integrating the second part of \eqref{graboundRem} along $\gamma$
and using \eqref{implelem}, we obtain
$$\begin{aligned}
&\log\Big(\sqrt{\big|\log\log\big(\ts\frac{N}{u(x_2,t)}\big)\big|}+K\Big)\le 
\log\Big(\sqrt{\big|\log\log\big(\ts\frac{N}{u(x_1,t)}\big)\big|}+K\Big)+C\ts\frac{d}{R}\\
&\ \ \Rightarrow \sqrt{\big(-\log\log\big(\ts\frac{N}{u(x_2,t)}\big)\big)}+K\le 
\lambda\Big(\sqrt{\big(-\log\log\big(\ts\frac{N}{u(x_1,t)}\big)\big)}+K\Big)\\
&\ \ \Rightarrow -\log\log\big(\ts\frac{N}{u(x_2,t)}\big)\le 
\lambda^4\big(-\log\log\big(\ts\frac{N}{u(x_1,t)}\big)\big)+(\lambda^2-1)K^2\\
&\ \ \Rightarrow \log\log\big(\ts\frac{N}{u(x_1,t)}\big)\le \lambda^{-4}
\log\log\big(\ts\frac{N}{u(x_2,t)}\big)+ (1-\lambda^{-4})K^2\\
&\ \ \Rightarrow \log\big(\ts\frac{N}{u(x_1,t)}\big)\le 
\big[\log\big(\ts\frac{N}{u(x_2,t)}\big)\big]^{\lambda^{-4}} e^{(1-\lambda^{-4})K^2}.
\end{aligned}$$
It follows that
$$\begin{aligned}
\big|\log\big(\ts\frac{u(y,t)}{u(x,t)}\big)\big|
&=\log\big(\ts\frac{u(x_2,t)}{u(x_1,t)}\big)
=\log\big(\ts\frac{N}{u(x_1,t)}\big)-\log\big(\ts\frac{N}{u(x_2,t)}\big)\\
&\le\big[\log\big(\ts\frac{N}{u(x_2,t)}\big)\big]^{\lambda^{-4}} e^{(1-\lambda^{-4})K^2}-\log\big(\ts\frac{N}{u(x_2,t)}\big)\\
&\le\big[\log\big(\ts\frac{N}{u(x,t)}\big)\big]^{\lambda^{-4}} e^{(1-\lambda^{-4})K^2}-\log\big(\ts\frac{N}{u(x,t)}\big),\end{aligned}$$
where we used \eqref{monotxi} and the fact $z\mapsto z^{\lambda^{-4}} e^{(1-\lambda^{-4})K^2}-z$
is nondecreasing for $z\in (0,1]$ due to $\lambda\ge 1$.
 Using $d\le R$ and then letting $\eps\to 0$ in $N=M+\eps$, we deduce 
\be{modulus2a}
\big|\log\big(\ts\frac{u(y,t)}{u(x,t)}\big)\big|
\le\big[\log\big(\ts\frac{M}{u(x,t)}\big)\big]^{\frac{1}{1+Cd/R}} e^{CdR(t^{-1}+k)}-\log\big(\ts\frac{M}{u(x,t)}\big),
\ee
hence
\be{modulus2}
\ts\frac{|u(y,t)-u(x,t)|}{u(x,t)}
\le\exp\Big\{\big[\log\big(\ts\frac{M}{u(x,t)}\big)\big]^{\frac{1}{1+Cd/R}} e^{CdR(t^{-1}+k)}-\log\big(\ts\frac{M}{u(x,t)}\big)
\Big\}-1.\ee
Now if we assume
$$d(x,y)\le \frac{cR}{R^2(t^{-1}+k)+|\log\log(\frac{M}{u(x,t)})|}$$
with $c=c(n)>0$ small, by a simple computation using $M/u(x,t)\le 2$, we in particular 
deduce from \eqref{modulus2a} that
 $\frac{u(y,t)}{u(x, t)}\ge \sqrt{2/e}$.
Consequently, $\rho=R/2$ in \eqref{defrhoH}
and this concludes the proof in the second case.
\qed

\subsection{Proof of Theorem~\ref{propBernAppl}}
Fix $t\in(0,T]$ and $x,y\in \overline{B_{R/2}}$ with $u(x,t)\le u(y,t)\le M$.
We may assume $u(x,t)<\min(u(y,t),M/2)$ since otherwise we are done.
Let $\gamma~:~[0,1]\to \mathcal{M}$ be a geodesic such that $\gamma(0)=x$ and $\gamma(1)=y$.
Let $\sigma_0=\min\{\sigma\in(0,1],\ \gamma(\sigma)=\min(u(y,t),M/2)\}$.
We have $\sigma_0\in(0,1]$ by continuity and $u(\gamma(\sigma),t)\le M/2$ for all $\sigma\in[0,\sigma_0]$. Set 
$z=\gamma(\sigma_0)$ and $\gamma_1=\gamma([0,\sigma_0])$.
Keeping the notation \eqref{notationKlambda} and applying \eqref{modulus1a} 
(and the subsequent sentence) with $x_1=z$ and $x_2=x$, we obtain
$$\log(M/u(x,t))\le \lambda^4 \log(M/u(z,t))+(\lambda^2-1)K^2.$$
Setting $\theta=\big(1+C\ts\frac{d}{R}\big)^{-1}$ and using $d\le R$, we deduce
$$ \log u(z,t)\le \theta \log u(x,t)+ (1-\theta)\log M+A,\quad A:=CR^2\big(t^{-1}+k\big),$$
hence
$$u(y,t)\le M\le 2\min(u(y,t),M/2) =2u(z,t)\le 2e^A M^{1-\theta}u^\theta(x,t),$$
and \eqref{pseudoHarnack} follows.
\qed

  \begin{remark} \label{remOpt1}
As mentioned after Theorem~\ref{propBernAppl}, the sharpness of \eqref{pseudoHarnack} 
 for $k=0$ can be seen on the family of solutions
$$u=u_\tau(x,t)=G(x_0,x,t+\tau),\quad (x,t)\in\mathcal{M}\times(0,\infty), \quad \tau>0,\ x_0\in \mathcal{M},$$
 where $G$ is the heat kernel of $\mathcal{M}$.
Let us first consider the case $\mathcal{M}={\R}^n$.
For each $R, \tau>0$, taking $x_0=0$ without loss of generality,
$Q=B_R\times(0,\tau]$, $y=0$, $|x|=d(x,y)=R/2$ and $t=\tau$,
we have 
  $u(x,t)=(2\tau)^{-n/2}e^{-R^2/32\tau}$, $u(y,t)=(2\tau)^{-n/2}$, $M=\tau^{-n/2}$, $\theta=\big(1+C)^{-1}$, 
$L=\exp\big(C\big(1+\frac{R^2}{\tau}\big)\big)$, with $C=C(n)>0$.
Therefore, \eqref{pseudoHarnack} yields the inequality
$$\begin{aligned}
1
&\le \frac{L M^{1-\theta}u^\theta(x,t)}{u(y,t)}\\
&=\ts\frac{2}{(2\tau)^{-n/2}}\exp\big(C\big(1+\frac{R^2}{\tau}\big)\big)\tau^{-n(1-\theta)/2}(2\tau)^{-n\theta/2}
\exp\big(-\frac{\theta R^2}{32\tau}\big)\\
&=2^{1+n(1-\theta)/2}e^C\exp\Big[\big(C-\ts\frac{1}{32(1+C)}\big)\ts\frac{R^2}{\tau}\Big],\quad\hbox{for all }R,\tau>0,
\end{aligned}$$
and the validity of the latter precisely depends on the magnitude of the constant $C=C(n)$:
 it is true if and only if
  $C\ge\ts\frac{3-2\sqrt 2}{4\sqrt{2}}$. 
 
 In the case of a general manifold with $k=0$, 
 the above argument applies with straightforward changes by using
the Li-Yau estimate  \eqref{group2} for the heat kernel.
  \end{remark}

\section{Proof of Proposition~\ref{prop-optim2}}\label{sect4}
In this proof we set $J=\sqrt{H}$ for convenience.
 Fix $s, t, R>0$ and set $Q=B_R\times (0,T]$ with $T:=t$.

$\bullet$ First consider the solutions $u=u_a(x,\tau)=e^{ax+a^2\tau}$ in $Q$.
Choosing $\tau=t$, $a=2s/R$ and $x=Re_1/2$, we obtain
$$M=e^{a+a^2T},\quad M/u(x,t)=e^{aR/2},\quad \log(M/u(x,t))=s,$$
   hence 
    \be{optim-case01}
    J(s,t,R)\ge \ts|\frac{\nabla u}{u}|=a=2R^{-1}s,
        \quad\hbox{ for $s, t, R>0$.}
    \ee
In particular, we also have
    \be{optim-case02}
    J(s,t,R)\ge 
    \begin{cases}
    t^{-1/2}s&\hbox{if $s\le 1$ and $t\ge R^2$}\\
    t^{-1/2}s^{1/2}&\hbox{if $s\ge 1$ and $s\ge R^2/t$.}
        \end{cases}
    \ee

$\bullet$ Now allowing $R\in(0,\infty]$, we next consider the solutions
$$u=u_a(x,\tau)=(\tau+a)^{-n/2}e^{-|x|^2/4(\tau+a)} \quad\hbox{in $Q$}.$$
We have 
$$M=a^{-n/2},\quad M/u=(\ts\frac{\tau+a}{a})^{n/2}e^{|x|^2/4(\tau+a)},\quad \frac{|\nabla u|}{u}=\frac{|x|}{t+a}.$$
If $1\le s\le R^2/t$, we choose $a=\tau=t$ and $x=\ts\frac12 \bar R e_1$,
with $\bar R=R$ if $R<\infty$ or $\bar R=\sqrt{st}$ otherwise.
We obtain
    \be{optim-case1}
    J(s,t,R)\ge \frac{|\nabla u|}{u}=\frac{|x|}{2t}=\frac{\bar R}{4t}\ge \frac14 \sqrt{\frac{s}{t}},
    \quad\hbox{ if $1\le s\le R^2/t$.}
    \ee
If $s\le 1$ and $t\le R^2$,
we choose $\tau=t$ and $x=\ts\frac12 \sqrt t e_1$ (hence $|x|\le R/2$).
Then there exists a unique $a>0$ such that
$$\psi(a):=\ts\frac{n}{2}\log(1+\ts\frac{t}{a})+\frac{t}{16(t+a)}\equiv \log(M/u(x,t))=s$$
(since $\psi$ is decreasing, with $\lim_{a\to 0}\psi(a)=\infty$ and $\lim_{a\to \infty}\psi(a)=0$).
In addition, we have $a\ge c(n)t$ and $s\le \ts\frac{n}{2}\ts\frac{t}{a}+\frac{t}{16(t+a)}\le C(n) \frac{t}{a}$. Consequently,
    \be{optim-case2}
    \frac{|\nabla u|}{u}\ge c(n)\frac{|x|}{a}=c(n)\frac{\sqrt{t}}{a}\ge c(n)\frac{s}{\sqrt{t}},
    \quad\hbox{ if $s\le 1$ and $t\le R^2$}.
    \ee
Gathering cases  \eqref{optim-case02}, \eqref{optim-case1} and \eqref{optim-case2}, we obtain $J(s,t,R)\ge ct^{-1/2}s^{1/2}$ for $s\le 1$ and  $J(s,t,R)\ge ct^{-1/2}s$ for $s\ge 1$.
Combining this with \eqref{optim-case01} and going back to definition \eqref{defH}, inequality \eqref{lowerg} follows.
As for \eqref{lowergB}, it is a consequence of \eqref{optim-case1} and \eqref{optim-case2} with $R=\infty$.
\qed


\medskip

\begin{thebibliography}{99}

\bibitem{Ai04}
S. Aida, 
{\it Precise Gaussian estimates of heat kernels on asymptotically flat Riemannian manifolds with poles},
Recent developments in stochastic analysis and related topics,
River Edge, NJ, World Scientific, 1-19 (2004).

\bibitem{ACMM}
D. Alonso-Or\'an, F. Chamizo, A.D. Mart\'inez, A. Mas, 
{\it Pointwise monotonicity of heat kernels}, 
Rev. Mat. Complutense 36 (2023),  207-220.

\bibitem{ATW}
M. Arnaudon, A. Thalmaier, F.-Y. Wang, 
{\it Gradient estimates and Harnack inequality on noncompact Riemannian manifolds},
Stochastic Processes Appl. 119 (2009), 3653-3670. 

\bibitem{Bai}
M. Bailesteanu, 
{\it Bounds on the heat kernel under the Ricci flow},
Proc. Am. Math. Soc. 140 (2012), 691-700.

\bibitem{BCP}
M. Bailesteanu, X. Cao, A. Pulemotov,
{\it Gradient estimates for the heat equation under the Ricci flow},
J. Funct. Anal. 258 (2010), 3517-3542.

\bibitem{BBG}
D. Bakry, F. Bolley, I. Gentil, 
{\it The Li-Yau inequality and applications under a curvature-dimension condition},
Ann. Inst. Fourier 67 (2017),  397-421.

\bibitem{bauer2015li}
F. Bauer, P. Horn, Y. Lin, G. Lippner, D. Mangoubi, S.-T. Yau,
{\it Li-Yau inequality on graphs}.
 J. Differential Geom. 99 (2015), 359-405.
 

 \bibitem{BM25}
A. Burchard, A.D. Mart\'\i nez, 
{\it On monotonicity of heat kernels: a new example and counterexamples},
Trans. Am. Math. Soc. 378 (2025), 6203-6225.


 \bibitem{BM26}
A. Burchard, A.D. Mart\'\i nez, 
{\it A new characterization of the round sphere},
Preprint, arXiv:2611.08751

  \bibitem{CCT21}
J. Cao, L.-J. Cheng, A. Thalmaier,
{\it Hessian heat kernel estimates and Calder\'on--Zyg\-mund inequalities on complete Riemannian manifolds},
Preprint, arXiv:2108.13058.

 \bibitem{CCM}
D. Castorina, G. Catino, C. Mantegazza,
{\it Semilinear Li and Yau inequalities},
Ann. Mat. Pura Appl. (4) 202 (2023), 827-850.


 \bibitem{CLW}
X. Chen, X.-M. Li, B. Wu,
{\it Logarithmic heat kernel estimates without curvature restrictions},
Ann. Probab. 51 (2023), 442-477.

 \bibitem{CY}
S.Y. Cheng, S.-T. Yau,
{\it Differential equations on
Riemannian manifolds and their geometric applications},
 Comm. Pure Appl. Math. 28 (1975), 333-354.

\bibitem{colding2021optimal}
T. Colding and W.P. Minicozzi II,
 {\it Optimal bounds for ancient caloric functions},
 Duke Math. J. 170 (2021), 4171-4182.
 

 \bibitem{Dri94}
B.K. Driver, 
 {\it A Cameron-Martin type quasi-invariance theorem for pinned Brownian motion on a compact Riemannian manifold},
Trans. Amer. Math. Soc. 342 (1994), 375-395.

\bibitem{DN}
H.T. Dung, T.D. Nguyen, 
{\it Sharp gradient estimates for a heat equation in Riemannian manifolds},
Proc. Am. Math. Soc. 147 (2019), 5329-5338.

\bibitem{Eng}
A. Engoulatov, 
{\it A universal bound on the gradient of logarithm of the heat kernel for manifolds with bounded Ricci curvature},
J. Funct. Anal. 238 (2006), 518-529.


\bibitem{EK22}
S. Eswarathasan, T. Kolokolnikov, 
{\it Laplace-Beltrami spectrum of ellipsoids that are close to spheres and analytic perturbation theory},
IMA J. Appl. Math. 87 (2022), 20-49.

\bibitem{Grigoryan91}
A.~Grigor'yan.
{\it The heat equation on noncompact Riemannian manifolds (Russian)},
 Mat. Sb. 182 (1991), 55-87 (English translation: Math. USSR Sb. 72 (1992), 47-77.

\bibitem{Ham93}
 R.S. Hamilton, 
 {\it A matrix Harnack estimate for the heat equation},
 Comm. Anal. Geom 1 (1993), 113-126. 
 
\bibitem{HZ}
Q. Han, Q.S. Zhang, 
{\it An upper bound for Hessian matrices of positive solutions of heat equations},
J. Geom. Anal. 26 (2016), 715-749.

\bibitem{Hsu}
E.P. Hsu, 
{\it Estimates of derivatives of the heat kernel on a compact Riemannian manifold},
Proc. Am. Math. Soc. 127 (1999), 3739-3744.

\bibitem{HZL}
G. Huang, Z. Huang, H. Li, 
{\it Gradient estimates for the porous medium equations on Riemannian manifolds},
J. Geom. Anal. 23 (2013), 1851-1875.

\bibitem{HuZ}
J.-G. Huang, H.-C. Zhang, 
{\it Localized elliptic gradient estimate for solutions of the heat equation on $RCD^*(K,N)$ metric measure spaces},
Manuscr. Math. 161 (2020),  303-324.

\bibitem{Kot}
B.L. Kotschwar, 
{\it Hamilton's gradient estimate for the heat kernel on complete manifolds},
Proc. Amer. Math. Soc. 135 (2007), 3013-3019.

\bibitem{Li91}
J.Y. Li, 
{\it Gradient estimates and Harnack inequalities for nonlinear parabolic and nonlinear elliptic equations on Riemannian manifolds},
J. Funct. Anal. 100 (1991), 233-256.

\bibitem{LX}
J. Li, X. Xu, 
{\it Differential Harnack inequalities on Riemannian manifolds I: Linear heat equation},
Adv. Math. 226 (2001), 4456-4491.

\bibitem{LZh}
L. Li, Z.L. Zhang, 
{\it On Li-Yau heat kernel estimate},
Acta Math. Sin., Engl. Ser. 37 (2021), 1205-1218.

\bibitem{Li12}
P.~Li,
{\it Geometric analysis}, volume 134 of  Cambridge Studies in Advanced Mathematics.
 Cambridge University Press, Cambridge, 2012.

\bibitem{LiYau86}
P.~Li and S.-T. Yau.
{\it On the parabolic kernel of the Schr\"odinger operator}.
 Acta Math.~156 (1986), 153-201.

\bibitem{LZ}
Y. Li, X.R. Zhu, 
{\it Li-Yau-Hamilton estimates and Bakry-Emery-Ricci},
J. Differ. Equ. 260 (2016), 3270-3301.

\bibitem{lin2019ancient}
F. Lin, Q.S. Zhang.
{\it On ancient solutions of the heat equation}, 
 Comm. Pure Appl. Math. 72 (2019), 2006-2028.

\bibitem{Liu}
S.P. Liu, 
{\it Gradient estimates for solutions of the heat equation under flow},
Pac. J. Math. 243 (2009), 165-179.

\bibitem{LNVV}
P. Lu, L. Ni, J.-L. V\'azquez, C. Villani, 
{\it Local Aronson-Benilan estimates and entropy formulae for porous medium
and fast diffusion equations on manifolds},
J. Math. Pures Appl. 91 (2009), 1-19.

\bibitem{Lu23}
Z. Lu,
{\it Differential Harnack inequalities for semilinear parabolic equations on Riemannian manifolds. 
I: Bakry-\'Emery curvature bounded below},
J. Differ. Equations 377 (2023), 469-518.

\bibitem{MZS}
L. Ma, L. Zhao, X.F. Song, 
{\it Gradient estimate for the degenerate parabolic equation $u_t = \Delta F(u) + H(u)$ on manifolds},
J. Differ. Equations 224 (2008), 1157-1177.

\bibitem{MaSt96}
P. Malliavin, D.W. Stroock,
{\it Short time behavior of the heat kernel and its logarithmic derivatives,}
J. Differential Geom. 44 (1996), 550-570.

\bibitem{mosconi2021}
S. Mosconi.
{\it Liouville theorems for ancient caloric functions via optimal growth conditions},
 Proc. Amer. Math. Soc., 149 (2021), 897-906.


\bibitem{NeSa}
R.W. Neel, L. Sacchelli,
{\it Localized bounds on log-derivatives of the heat kernel on incomplete Riemannian manifolds},
Ann. Inst. Henri Poincaré, Probab. Stat. 62 (2026), 505-522.


\bibitem{LT}
L. Ni and L. F. Tam, 
{\it K\"ahler-Ricci flow and the Poincar\'e-Lelong equation},
Comm. Anal. Geom 12 (2004), 111-141.


\bibitem{LSC02}
L. Saloff-Coste, 
{\it Aspects of Sobolev-type inequalities},
London Mathematical Society Lecture Note Series. 289. Cambridge University Press, 2002.

\bibitem{souplet2006sharp}
Ph. Souplet, Q.S. Zhang,
{\it Sharp gradient estimate and {Y}au's {L}iouville theorem for the heat equation on noncompact manifolds}.
 Bull. London Math. Soc. 38 (2006), 1045-1053.

\bibitem{STu}
D. Stroock, J. Turetsky, 
{\it Upper bounds on derivatives of the logarithm of the heat kernel},
Commun. Anal. Geom. 6 (1998), 669-685.


\bibitem{Uhl}
 Uhlenbeck, K., 
{\it Generic properties of eigenfunctions},
Amer. J. Math. 98 (1976), 1059-1078.

\bibitem{Wa}
L.F. Wang, 
{\it Gradient estimates on the weighted $p$-Laplace heat equation},
J. Differ. Equation 264 (2018) 506-524.

\bibitem{WWa}
W. Wang,
{\it Heat kernel under integral Ricci curvature condition and Riesz transform},
J. Differ. Equations 440 (2025), 113432.

\bibitem{Wu}
J.Y. Wu, 
{\it Elliptic gradient estimates for a nonlinear heat equation and applications},
Nonlinear Anal. 151 (2017), 1-17.
 
\bibitem{Ya}
Y.Y. Yang, 
{\it Gradient estimate for a nonlinear parabolic equation on Riemannian manifold},
Proc. Am. Math. Soc. 136 (2008) 4095-4102.

\bibitem{ZhIMRN}
Q.S. Zhang, 
{\it Some gradient estimates for the heat equation on domains and for an equation by Perelman},
Int. Math. Res. Not. 2006, Article ID 92314, 39 p. 

\bibitem{ZhSurvey}
Q.S. Zhang, 
{\it Log gradient estimates for heat type equations on manifolds},
Preprint arXiv:2407.20719 (2024).

\bibitem{ZZ}
Q.S. Zhang, M. Zhu, 
{\it Li-Yau gradient bounds on compact manifolds under nearly optimal curvature conditions},
J. Funct. Anal. 275 (2) (2018), 478-515.


\bibitem{Zh}
X.B. Zhu, 
{\it Gradient estimates and Liouville theorems for nonlinear parabolic equations on noncompact Riemannian manifolds},
Nonlinear Anal. 74 (2011), 5141-5146.

\end{thebibliography}
\enddocument